\makeatletter
\def\@maketitle{%
\defaultfont\normalsize
\let\@makefnmark\relax \let\@thefnmark\relax \ifx\@empty\@subjclass\else
\@footnotetext{1991 {\it Mathematics Subject
Classification}.\enspace
\@subjclass.}\fi
\ifx\@empty\@keywords\else
\@footnotetext{{\it Key words and phrases.}\enspace \@keywords.}\fi
\ifx\@empty\@thanks\else
\@footnotetext{\@thanks}\fi
\topskip66\p@ 
\vtop{\centering{\baselineskip14\p@\bf
\expandafter{\@title}\@@par}%
\global\dimen@i\prevdepth}%
\prevdepth\dimen@i
\ifx\@empty\@authors
\else
\baselineskip32\p@
\vtop{\@andify{ AND }\@authors
\centering{{\@authors}\@@par}%
\global\dimen@i\prevdepth}\relax
\prevdepth\dimen@i
\fi
\ifx\@empty\@dedicatory
\else
\baselineskip18\p@
\vtop{\centering{\small\it\@dedicatory\@@par}%
\global\dimen@i\prevdepth}\prevdepth\dimen@i \fi
\ifx\@empty\@date\else
\baselineskip24\p@
\vtop{\centering\@date\@@par
\global\dimen@i\prevdepth}\prevdepth\dimen@i \fi
\normalsize
\dimen@32\p@ \advance\dimen@-\baselineskip \vskip\dimen@\@plus14\p@
} 
\makeatother

\documentstyle[amscd,amssymb,verbatim,12pt,leqno]{amsart}

\theoremstyle{plain}
\newtheorem{Thm}{Theorem}[section]
\newtheorem{Cor}[Thm]{Corollary}
\newtheorem{Lem}[Thm]{Lemma}
\newtheorem{Prop}[Thm]{Proposition}
\newtheorem{Proc}[Thm]{Procedure}
\newtheorem{Prob}[Thm]{Problem}
\newtheorem{Main}{Theorem}
\renewcommand{\theMain}{\kern -6pt}

\theoremstyle{definition}
\newtheorem{Ex}[Thm]{Example}
\newtheorem{Rem}[Thm]{Remark}

\newtheorem{Def}[Thm]{Definition}

\newcommand{\ra}{\rightarrow}
\newcommand{\Spec}{\operatorname{Spec}}

\renewcommand{\O}{{\Cal O}}

\newcommand{\Hom}{\operatorname{Hom}}

\def\morph#1{\overset{#1}{\ra}}

\def\dsp#1{$\displaystyle{#1}$}
\def\Cal#1{{\cal #1}}
\def\BBB#1{{\Bbb #1}}
\def\d{{\Cal D}}
\def\Frak #1{{\frak{#1}}}
\def\pB{\B_W((P))}
\def\pBf{\pB_{fin}}
\def\mat#1{M_{{#1}}(\B_W)}
\def\matp#1{M_{{#1}}(\pBf)}
\def\qc{quasicoherent}
\def\gp{Y}
\def\ft{\Phi}
\def\lie{\Frak g}
\def\E{\Cal E}
\def\map{\phi}

\theoremstyle{remark}
\newtheorem{notation}{Notation}[section] 

\numberwithin{equation}{section}

\newcommand{\Oplus}{\operatornamewithlimits{\oplus}}

\newcommand{\A}{{\Cal A}}
\newcommand{\B}{{\Cal B}}
\renewcommand{\O}{{\Cal O}}
\newcommand{\F}{{\Cal F}}
\newcommand{\G}{{\Cal G}}

\newcommand{\D}{{\bold D}}

\newcommand{\Pic}{\operatorname{Pic}}
\newcommand{\Sym}{\operatorname{Sym}}

\newcommand{\End}{\operatorname{End}}
\newcommand{\Ext}{\operatorname{Ext}}

\newcommand{\gr}{\operatorname{gr}}

\renewcommand{\P}{\Cal P}

\newcommand{\C}{\BBB C}

\title[Dynamics of the Krichever construction in several variables]{Dynamics of 
the Krichever construction in several variables}

  \author{ Mitchell Rothstein}

\begin{document}

\begin{abstract}The  Krichever construction in one variable,  that is, for spectral curves,  
linearizes the KdV-hierarchy on the jacobian of the curve.
We carry out an appropriate generalization of the Krichever construction for an arbitrary 
projective variety $X$ and  determine the corresponding
nonlinear dynamics, which are then  linearized  on an extension of $\Pic^0(X)$.
\end{abstract}

\maketitle

\section{Introduction}\label{sec:intro}    \subsection{}The 
algebro-geometric inverse scattering method 
for solving nonlinear partial differential equations is well-developed in the 
case of one-dimensional
spectral varieties.   Given a nonlinear 
evolution equation,  
one interprets the evolving  fields
as operators acting on objects naturally associated with the geometry of 
the curve. The nonlinear dynamics  then becomes
linear flow on the Picard variety of the  curve.
The best-known example is the 
{\em Krichever construction},   
which 
associates to a projective curve, a smooth point and  some
additional data,  
a solution of the KdV-hierarchy, \cite{K1,K2, K3, KN,  Mum, S,SW}.

Once Krichever's construction appeared,  it was natural to ask  whether a similar construction could be used to linearize
dynamics on Picard varieties of higher dimensional varieties,  \cite{S,T}.   An important,  
though nondynamical contribution was made by Nakayshiki, \cite{N},
who used the Fourier-Mukai transform, \cite{Muk}, to  represent the ring of 
functions on 
a principally polarized abelian variety $X$, with poles along the 
theta divisor,   as matrix differential operators in $\dim(X)$-many variables.   
This example is nondynamical in the sense that $\dim(X)$ is also
$\dim(\Pic(X))$,   (which is just $X$ in this case.)    In the case of a curve 
of genus $g$,  
on the other hand,  one has $g-1$ independent KdV-flows.   In
\cite{R}, a
 general setting was described,  in which one may recover not only 
Nakayshiki's example, 
but other examples as well.  For instance, one could replace the
curve by the family of lines on a smooth cubic threefold,  
(this is the Fano surface,) and 
replace the point by an incidence divisor.   This example is
interesting from the standpoint of dynamics,  because  
$\dim(\Pic(X))-\dim(X)=5-2=3$,  so one anticipates that there
will be three commuting flows, although the flows were not given 
in that paper.

As it turns out, the situation is quite a bit more general than it
first seemed.   In the present work, we will carry out the Krichever 
construction with  any Cohen-Macauley projective variety $X$,  and 
determine the
corresponding dynamics.    The flows will be defined on an extension 
$Y\to\Pic^0(X)$, of $\Pic^0(X)$ by several copies of $\BBB G_m$.

 \subsection{} To  describe the flows,  let $n=\dim(X)$, and introduce the following additional 
data: 

$\bullet$ $Z\subset X$, a reduced, ample Cartier hypersurface.   Set $\Cal N=\O(Z)/\O$, the normal bundle to $Z$, regarded as a locally
free sheaf on $Z$.

$\bullet$ an $n$-dimensional basepoint free subspace $W_0\subset H^0(Z,\Cal N)$.    Let $\phi:Z\to\BBB P(W_0^*)$ denote the
corresponding morphism.

$\bullet$  a point $P\in\BBB P(W_0^*)$, such that $\phi$ is \'etale in a 
neighborhood of $\phi^{-1}(P)$.    Set $S=\phi^{-1}(P)$ and let $\Cal I_S$ denote the ideal sheaf of $S$ in $X$.

$\bullet$   a choice of $n$-dimensional subspace $W\subset H^0(X,\O(Z)/\Cal I_S)$,  
mapping isomorphically onto $W_0$.

$\bullet$ a basis $\xi_1,...,\xi_n$ for $W$, such that the images of $\xi_2,...,\xi_n$ in $W_0$ vanish  on $S$.

Our group $Y$ is then the group of degree-zero line bundles on $X$,  
equipped with a trivialization along $S$.   There is a natural map
\begin{equation}\label{eq:natural map}
W\to H^1(X,\Cal I_S)=\text{invariant
vector fields on
$Y$.}\end{equation}
Denoting by $\partial_i$ the image of $\xi_i$ under the map \eqref{eq:natural map}, one has an algebra 
$\O_Y[\xi_1,...,\xi_n]$, in which the $\xi $'s commute among themselves, and satisfy $[\xi_i,f]=\partial_i(f)$ for $f\in\O_Y$.

In section \ref{sec:fmt} we introduce a functor $\ft $, a version
of the Fourier-Mukai transform, that carries $\O_X(*Z)$-modules to $ \O_Y[\xi_1,...,\xi_n]$-modules.

\begin{Main}[Theorem \ref{Locally Free at the generic point}  and Corollary \ref{cor:CM}] Assume $X$ is Cohen-Macauley.  
Let $\Cal M$,  $\Cal L$ and $\Cal F'$ be line bundles,  with $\Cal M$ and $\Cal L$ ample.
Then there exist 
integers
$a$ and
$b$, such that for every collection of data as described above,  if 
$\O(Z)\simeq\Cal M^b$ and if we set $\F=\Cal L^a\otimes\F'$,  there exists a Zariski-dense open subset
$\Cal U\subset Y$, such that 
 the $ \O_Y[\xi_1,...,\xi_n]$-module $\ft(\Cal F(*Z))$ is free on $\Cal U$.    Let $d$ denote the degree
of
$Z$.   Then the rank of $\ft_S(\Cal F(*Z))$ is $d$. \end{Main}

  By functoriality,   one has a map, which turns out to be an imbedding,
\begin{equation}\label{eq:abstract imbedding}H^0(X,\O(*Z))\to \End_{\O_Y[\xi_1,...,\xi_n]}(\ft(\Cal F(*Z)))\ .\end{equation}
A choice of basis for $\ft(\Cal F(*Z))|_{\Cal U}$ then gives
an imbedding
\begin{align}\label{eq:imbedding}H^0(X,\O(*Z))&\to H^0(\Cal U,M_d(\O_Y[\xi_1,...,\xi_n]))\\
f&\mapsto L_f\ ,\notag\end{align}
into the  ring of $d\times d$
matrices over 
$H^0(\Cal U,\O_Y[\xi_1,...,\xi_n])$.   
This generalizes the Burchnall-Chaundy imbedding  of the ring
of functions on a curve with poles at point into the ring of differential operators in one variable,  \cite{BC}.   

The dynamics problem is to find,  for every 
translation invariant vector field  $\delta$  on $Y$, an explicit formula for $\delta(L_f)$.  
Since $L_f$ is
determined only up to conjugation, any such formula will depend on choices.  When $X$ is a curve and $Z$ is a point, there is a natural 
normalization that leads to the KdV-hierarchy, \cite{Mum}.
We shall show that in the general case, there are
many ways to make such choices and still obtain explicit formulas for the flow. These choices will also
affect the number of copies of
$\BBB G_m$ entering into  the extension $Y$, since it may be possible to throw away the trivializations at some points of $S$.   In some cases it will be possible
to take
$Y=\Pic^0(X)$,  which is of particular interest when
$\Pic^0(X)$ is an abelian variety.
   The flow that can be described in the fewest words uses all of $S$.   It is discussed in section \ref{sec:flows}.
The more exotic flows that use smaller extensions are discussed in section \ref{sec:more flows}. 

The flow given in section \ref{sec:flows} goes as follows.  Invert $\xi_1$ by allowing formal power series in $\xi_1^{-1}$, 
(microlocalization):
$$H^0(\Cal U,M_d(\O_Y[\xi_1,...,\xi_n]))\subset H^0(\Cal U,M_d(\O_Y[\xi_2,...,\xi_n]((\xi_1^{-1})))) $$   
For all nonzero $f\in H^0(X,\O(*Z))$,  $L_f$ is invertible in this larger ring of matrices.   

\begin{Main}[Theorem \ref{draconian}] For an appropriate choice of   basis for $\ft_S(\Cal F(*Z))|_{\Cal U}$,  the imbedding
\eqref{eq:imbedding} is such that for every translation invariant vector field  
$\delta$, there exist $f,g\in H^0(X,\O(*Z))$  such that for all $h\in H^0(X,\O(*Z))$, 
\begin{equation}\label{eq:flows zero}
\delta(L_h)=[({L_g}^{-1}{L}_f)_+,L_h]\ .\end{equation} \end{Main}

For the modified flows discussed in section \ref{sec:more flows},  the diagonal entries of  
$({L_g}^{-1}{L}_f)_+$  are altered by  adding  terms  from the negative part of ${L_g}^{-1}{L}_f$.

When  $n=d=1$,  we recover the KdV-hierarchy as follows.  The variety $X$ is now a curve,  and $Z$ is a point.  Choose a local coordinate $z$ at $Z$, 
in such a way that for some integer $r$ there exists $h\in H^0(X,\O(*Z))$ such that $z^r=h$.  For each integer $j$,   $1/z^j$ defines a vector
field
$\delta_j$ on $\Pic^0(X)$. Then,  taking $\delta=\delta_j$ in equation \eqref{eq:flows zero},  $f$ and $g$ are such that
$f/g=1/z^{j}+o(z) $.  Then $({L_g}^{-1}{L}_f)_+=(L_h^{j/r})_+$. Then \eqref{eq:flows zero} is the KdV-hierarchy.

For general $n$ and $d$, it is not possible to recover $({L_g}^{-1}{L}_f)_+$ in the form $(L_h^{j/r})_+$.    Instead, for each vector field $\delta$, we have an
autonomous nonlinear system 
\begin{align}\label{eq:flows half}
&\delta(L_f)=[({L_g}^{-1}{L}_f)_+,L_f]\ .\notag\\&\delta(L_g)=[({L_g}^{-1}{L}_f)_+,L_g]\\
&[L_f,L_g]=0\ .\notag
\end{align}
Then, given another vector field $\delta_1$,  we get {commuting} flows for
 four commuting matrices, and so on.

On the other hand, in any particular example,  one may choose a set of generators $f_1,...,f_m$ for  $H^0(X,\O(*Z))$ as an algebra over $\BBB C$.   Then,  as
$\delta$ varies, one may view
\eqref{eq:flows zero} as a collection of commuting flows for the $m$ commuting matrices $L_{f_i}$,  with $L_f$ and $L_g$ now being  polynomials in 
$L_{f_1},...,L_{f_m} $. 

\subsection{} An important aspect of the one-variable theory is the {\em Krichever map}, $\{$Krichever data$\}\to\{$point in infinite
dimensional grassmannian$\}$,   discussed by Sato, \cite{S} and Segal and Wilson, \cite{SW}.  Progress has been made on
understanding the Krichever map  in several variables,  by Parshin, \cite{P}, for surfaces and by Osipov, \cite{O}, in higher dimensions. 
I will not  attempt in this paper to obtain the dynamics \eqref{eq:flows zero}
from flows on an infinite dimensional manifold,   thought it does seem that  such a construction might be possible.  
Other open problems include finding soliton solutions to  \eqref{eq:flows
zero} and finding explicit formulas in terms of theta functions.

\subsection{} The outline of the paper is as follows. Section \ref{sec:fmt} 
gives the basic 
results connecting the Fourier-Mukai transform  with the theory of 
$\Cal D$-modules.   Section
\ref{sec:cmdo} applies this theory, to obtain commuting matrix 
partial differential
operators.   Section \ref{sec:good basis} addresses the problem of 
choosing a good basis.  
Section \ref{sec:microdiff} brings in  microdifferential operators.  
The dynamics  are  obtained in section \ref{sec:flows}.  
In section
\ref{sec:more flows} we obtain the modified dynamics that descend to 
quotients of $Y$, gotten by  
throwing away the trivializations at some of the points of
$\map^{-1}(P)$.  Examples in which  the flow descends to $\Pic^0(X)$ 
itself
are given in section
\ref{sec:examples}.

\bigskip
It is a pleasure to acknowledge invaluable discussions with  Alexander 
Polishchuk and Robert Varley.

\section{Fourier-Mukai transform and $\Cal D$-modules}\label{sec:fmt}

\subsection{Fourier-Mukai transform}\label{fmt}

Let $X$ be a projective variety over $\BBB C$,   $S\subset X$ a 
nonempty finite subset, and $\gp$ the group of degree-zero line 
bundles on $X$, trivialized
along
$S$.     For a sheaf of 
$\O_X$-modules $\F$,  denote by $\F|_S$ the pullback of $\F$ 
to $S$. 
Let $\P$ denote the Poincar\'e line bundle on $X\times \gp$,   i.e.
a universal line bundle  equipped with a trivialization
$\Cal P|_{S\times
\gp}\simeq\O_{S\times
\gp}$.   One has a canonical isomorphism
\begin{equation}\label{square}
(id\times\pi_1)^*(\P)\otimes(id\times\pi_2)^*(\P)\simeq
(id\times m)^*(\P)\ ,
\end{equation}
where $\pi_1,\pi_2$ and $m$ are the projection and
multiplication maps, respectively,   on $\gp\times \gp$.  The Fourier-Mukai transform is defined as follows.

\begin{Def}\label{transform}\cite{Muk} 
Given an object 
$\F$ in the derived category of
\qc\  
$\O_X$-modules,  set
\begin{equation}
R\ft(\F)=R\pi_{{\gp}*}(\pi_X^*(\F)\otimes\P)\
,\end{equation}
where $\pi_X$,  $\pi_{\gp}$ are the two projection maps on $X\times 
\gp$. 
Set
\begin{equation}
\ft(\F)=H^0(R\ft(\F))\ .
\end{equation}
 \end{Def}

Let $\d$ denote the sheaf of linear differential operators on $\gp$.
One might  ask how to augment the data of an $\O_X$-module structure 
on a sheaf $\F$ so that $\ft(\F)$ will inherit a
$\d$-module structure (cf. \cite{L,R}.)

\begin{Def} An extension rigidified along $S$ is a locally split 
exact sequence of $\O$-modules
\begin{equation}
0\to\Cal F\to \Cal G\to \Cal H\to 0
\end{equation}
equipped with an isomorphism of extensions
\begin{equation}
\Cal F|_S\oplus\Cal H|_S  \simeq \Cal G|_S\ .
\end{equation}
\end{Def}
Let $\Cal I_S$ denote the ideal sheaf of $S$.
The extensions of $\O$ by $\O$ rigidified along $S$ are classified by 
$H^1(X,\Cal I_S)$.  Thus  one may make the following definition.

\begin{Def}    Let $\lie=H^1(X,\Cal I_S)$.  Let $\E$ denote 
the universal extension
\begin{equation}\label{extension}
0\to\Cal O\morph{t} \E\morph{p} \lie\otimes_{\C}\Cal O\to 0
\end{equation}
rigidified along $S$.
  Set
\begin{equation}
\A=\Sym_{\O}(\E)/\langle \
t-1 \rangle\ .
\end{equation}

\end{Def}

\begin{Prop}
$\ft(\cdot )$ defines a functor from   $\A$-modules to  $\d$-modules.
\end{Prop}

\begin{pf} The construction generalizes the one
given in \cite{R}.
Let
$\lie'=\Sym\lie^*/\langle\Sym^2\lie^*\rangle=\C\oplus\lie^*$.  Referring to   
\eqref{extension},
the map $\lie^*\to \Hom(\E)$, $\omega\mapsto 
t(\langle\omega,p(\cdot)\rangle)$, defines an action of   $\lie'$ on $\E$.
Thus $\E$  is the direct image of a line bundle on $X\times 
\Spec(\lie')$.  Denote this line bundle by  $\Cal P'$.
On the
other hand,
$\Spec(\lie')$ is the  first order neighborhood of the identity 
in $\gp$.   Then   $\Cal P'=\Cal P|_{X\times \Spec(\lie')}$
canonically.

It is  easy to see
that the following data are
equivalent:
\bigskip

\noindent DATA I.\ \    A splitting of the sequence
\begin{equation}\label{Fsequence} 0\to\F\to
\E\otimes_{\O}\F\to
\F\otimes_{\BBB C}\lie\to 0\ ,
\end{equation}
\noindent DATA II.\ \
An isomorphism
\begin{equation}
p^*(\F)\morph{\sim} p^*(\F)\otimes_{\O_{X\times \Spec(\lie')}}\P'
\end{equation}
restricting to the identity on $X$,   where $p$ is
projection onto $X$.

The first order neighborhood of the diagonal in $\gp\times \gp$ is 
canonically isomorphic to $\gp\times \Spec(\lie')$.   Denote by 
$p_i$,
$i=1,2$, the  two projections $\gp\times \Spec(\lie')\to \gp$ 
coming from this identification.   There is a third map, 
$\epsilon:\gp\times
\Spec(\lie')\to \gp$, given by projection onto $\Spec(\lie')$ 
followed by the inclusion map.

By equation \eqref{square},  one has a canonical isomorphism
\begin{equation}\label{littlesquare}
(id\times p_1)^*(\P)\simeq (id\times p_2)^*(\P)\otimes
(id\times \epsilon)^*(\P)\ .\kern -10pt
\end{equation}

It then follows from the base change formula that an instance of DATA 
II gives an isomorphism
\begin{equation}\label{connection}
p_1^*(\ft(\F))\simeq p_2^*(\ft(\F))
\end{equation}
restricting to the identity along ${\gp}$, that is to say, a 
connection on $\ft(\F)$.

For any splitting \dsp{\F\otimes_{\BBB C}\lie\morph{\tau}\E\otimes\F} of
\eqref{Fsequence} one may define its  curvature
  as follows.   Given an open set $U\subset X$, any splitting
\dsp{\O_U\otimes_{\BBB C}\lie\morph{\sigma}\E|_U} defines 
a splitting, $1\otimes\sigma$, of \eqref{Fsequence} restricted to 
$U$.   The
difference   is a section
\begin{equation}\tau-1\otimes\sigma\in
\lie^*\otimes\Gamma(U,\Cal End(\F))\ .
\end{equation}
The commutator is independent of the choice of $\sigma$,   so one has 
a  globally defined curvature
\begin{equation}[\tau,\tau]\in
\Lambda^2(\lie^*)\otimes\End(\F)\ .
\end{equation}
Then  the following data are equivalent:
\bigskip

splitting of \eqref{Fsequence} with vanishing curvature$\
\
\Leftrightarrow\ \ \A
$-module structure on $\F$.
\bigskip

It is easy to check that  $\ft$ carries the
curvature of the splitting to the curvature of the
connection.   Thus an $\A
$-module structure on $\F$ endows $\ft(\F)$ with a connection 
whose curvature vanishes,   i.e. a $\d$-module structure.
\end{pf}

\subsection{$\Cal D$-module structure on $\ft(\O_S)$}\label{alongS}    {\em A priori,} $\ft(\O_S)$ 
has two $\d$-module structures.   On the one hand, 
the pullback of $\A$ to $S$ is canonically trivial,
\begin{equation}\label{restrict}
\A|_S\simeq\O_S\otimes_{\C}\Sym_{\C}{\lie}\ .\end{equation}
Then there is a homomorphism
$\A\to\O_S$, sending
\hbox{$\lie\to 0$. } Then $\O_S$ is an $\A$-module so $ 
\ft(\O_S)$ becomes a $\d$-module.   On the other hand,
$\ft(\O_S)$ is canonically trivial,
\begin{equation}
\ft(\O_S)=H^0(S,\O_S)\otimes_{\C}\O_{\gp}\
\end{equation}
so it has its standard   $\d$-module structure.    We leave it 
to the reader to verify the following result.

\begin{Lem}\label{standard}   The action of $\d$ on 
$\ft(\O_S)$ induced by the canonical homomorphism  $\A\to\O_S$ 
is the
standard action of $\d$ on the trivial bundle.
\end{Lem}

\subsection{Algebras $\A$ and $\B$ associated to a linear map $M\to\lie $}\label{subsect:generalize}   The construction in   
subsection  \ref{fmt} may be 
generalized as follows.   Let $M$ be a $\C$-vector space and let
$M\morph{\tau}\lie $ be a linear map.   Then the extension $\E$ may 
be composed with $\tau$ to yield an extension

\begin{equation}\label{Mextension}
0\to\Cal O\morph{t} \Cal E_{\tau}\to M\otimes_{\C}\Cal O\to 0\ .
\end{equation}
On the other hand, $\tau$ may be viewed as a map from the trivial bundle  
$M\otimes_{\C}\Cal O_Y$ to the tangent sheaf of $Y$.  Thus  
$M\otimes_{\C}\Cal O_Y$ inherits  the structure of a Lie algebroid.
(See \cite{PR}.)

One then constructs sheaves of algebras $\A_{\tau}$ and $\B_{\tau}$ 
over $\O_X$ and $\O_{\gp}$ respectively,  in complete analogy with 
$\A$ and
$\d$:

\begin{Def} \ \

  Let
\dsp{
\A_{\tau}=\Sym_{\O}(\Cal E_{\tau})/\langle \
t-1 \rangle}.

For $m\in M$,   let $\partial_m$
denote the translation invariant  vector field on $\gp$  associated to $\tau(m)$.  
Let
$\B_{\tau}$ denote the universal enveloping algebra of the Lie algebroid
$M\otimes_{\C}\Cal O_Y$.   Thus, $\B_{\tau}$ is the
sheaf of associative
$\O_{\gp}$-algebras  generated by $M$,  with relations
\begin{align} [ m,f]&=\partial_m(f)\ ,\ \ f\in\O_{\gp},\ m\in M,\\
[{m_1},{m_2}]&=0\ ,\ m_1,m_2\in M\ .\end{align}
\end{Def}

\begin{Rem}If $\tau$ is the $0$-map,  then
\begin{align}\A_{\tau}&=\O_X\otimes\Sym M\\
\B_{\tau}&=\O_Y\otimes\Sym M\ \end{align}
as algebras.
\end{Rem}

\begin{Rem}\label{rightinverse}   There are canonical homomorphisms
\begin{equation}\A_{\tau}\to \A\ \ ,\ \
\B_{\tau}
\to \d\ .\end{equation}
On the other hand,   if $\tau$ is surjective,  then every right 
inverse of $\tau$ defines homomorphisms
\begin{equation}\A\to \A_{\tau}\ \ ,\ \
\d
\to \B_{\tau}\ .\end{equation}

\end{Rem}

\begin{Prop}
$\ft(\cdot)$ defines a functor from   $\A_{\tau}$-modules to 
$\B_{\tau}$-modules.
\end{Prop}

\begin{notation}  Every $\B_{\tau}$-module structure on an 
$\O_{\gp}$-module $\G$ determines,  and is determined by, a 
$\C$-linear
map $M\morph{\nabla}\End_{\C}(\G) $ such that for  all $f\in\O_{\gp}$ 
and all $m,m'\in M$, $ [\nabla_m,f]=\partial_m(f)$ and
\hbox{$[\nabla_{m},\nabla_{m'}]=0$.}   We will often refer to 
$\nabla$  itself as a $\B_{\tau}$-module structure or an action of 
$\B_{\tau}$.  \end{notation}

\subsection{Difference of two $\A_{\tau}$-module structures}   Let $\F$ be a sheaf of $\O_X$-modules equipped with 
two $\A_{\tau}$-module structures,  $\alpha_1$ and $\alpha_2$. 
Their
difference is a map
\begin{equation}
\alpha_1-\alpha_2:M\to\End(\F)\ .
\end{equation}
This in turn gives us a map 
\dsp{\ft(\alpha_1-\alpha_2):M\to\End(\ft(\F))}.

The proof of the following lemma is left to the reader.

\begin{Lem}\label{easylemma}\

 Let $\nabla^1$ and $\nabla^2$ be the  $\B_{\tau}$-module 
structures on  $\ft(\F)$  associated to $\alpha_1$ and $\alpha_2$ 
respectively.
Then
\begin{equation}
\nabla^1-\nabla^2=-\ft(\alpha_1-\alpha_2)\ .
\end{equation}

\end{Lem}

\begin{Rem}   The minus sign in part 1 of the lemma  is  a matter of 
convention, introduced for later convenience.
\end{Rem}

\subsection{\v Cech description}  Let $U_1,...,U_a$ be an affine open cover of $X$, and let $x(1),...,x(b)$ be a basis for $\lie$.  Represent each $x(i)$
as a cocycle $\{c_{\alpha,\beta}(i)\in H^0(U_{\alpha}\cap U_{\beta},\Cal I_S)\} $.  Then there exists sections $x_{\alpha}(1),...,x_{\alpha}(b)\in
H^0(U_{\alpha},\A)
$, such that $\A|_{U_{\alpha}}$ is freely generated as an algebra  over $\O|_{U_{\alpha}}$ by $x_{\alpha}(1),...,x_{\alpha}(b)
$. Furthermore,
\begin{equation}\label{eq:cocycle}
x_{\alpha}(i)|_{U_{\alpha}\cap U_{\beta}}- x_{\beta}(i)|_{U_{\alpha}\cap U_{\beta}}=c_{\alpha,\beta}(i)\ .\end{equation}
In particular,  $\O|_{U_{\alpha}}$  becomes an $\A$-module,  by letting $x_{\alpha}(i)$ act by $0$ for all $i$.  Then $\ft(\O|_{U_{\alpha}})$ is a $\Cal
D$-module.
 By functoriality, $H^0(U_{\alpha},\O) $ acts on $\ft(\O|_{U_{\alpha}})$, $\Cal D$-linearly.  
Therefore,  the
$\Cal D$-module structure on $\ft(\O|_{U_{\alpha}})$ may be identified canonically with a relative connection $\nabla^{\alpha} $ on $\P|_{U_{\alpha}\times
Y}$,  with 
$U_{\alpha}$ as the base. Moreover, $\O|_{U_{\alpha}\cap U_{\beta}} $ is an $\A$-module in  two ways, the difference between them being the map
\begin{align}\lie&\to H^0(U_{\alpha}\cap U_{\beta},\Cal I_S)\notag\\x(i)&\mapsto c_{\alpha,\beta}(i)\ .\end{align}  Then by lemma \ref{easylemma},
\begin{equation}\label{eq:relative cxn}
\nabla^{\alpha}_{\partial_i}-\nabla^{\beta}_{\partial_i}=-c_{\alpha,\beta}(i)\ ,\end{equation}
where $\partial_i$ is the vector field associated to $x(i)$.   Note that if we restrict $\nabla^{\alpha} $ to $\P|_{(U_{\alpha}\cap S)\times
Y}$, we get the standard connection on the trivial bundle.

Now let $\F$ be an $\A$-module.   Then $\F|_{U_{\alpha}}$ is an $\A$-module in two ways, the latter being the one where $x_{\alpha}(i)$ acts by $0$.
Taking the difference between these $\A$-module  structures, one finds that  the $\A$-module structure on $\F$  is described by a collection of commuting
endomorphisms $$\{\phi_{\alpha}(i)\in H^0(U_{\alpha},\End(\F))\} $$ such that
\begin{equation}\label{eq:coboundary}
\phi_{\alpha}(i)-\phi_{\beta}(i)=\text{multiplication by\ }c_{\alpha,\beta}(i)\ .\end{equation}
The relative connection $\nabla^{\alpha} $ on $\P|_{U_{\alpha}\times Y}$ induces a 
relative connection $\nabla^{\alpha} $ on $\pi_X^*(\F)\otimes\P|_{U_{\alpha}\times Y}$.    The formula
\begin{equation}\label{eq:local  formula}
\nabla_{\partial_i}=\nabla^{\alpha}_{\partial_i}+\phi_{\alpha}(i)\end{equation}
then defines a global relative connection on $\pi_X^*(\F)\otimes\P$,  with $X$ as the base.
The $\Cal D$-module structure on $\ft(\F)$ is then the direct image of this relative $\Cal D$-module structure.

If $M\morph{\tau}\lie$ is a linear map,  the situation is the same.
$$\A_{\tau}|_{U_{\alpha}}\simeq \O|_{U_{\alpha}}\otimes \Sym_{\C}(M_{\alpha})\ , $$
where $M_{\alpha}$ is a copy of $M$, such that
\begin{equation}\label{eq:Mcoboundary}
m_{\alpha}-m_{\beta}=c_{\alpha,\beta}(m)\ ,\end{equation}
where $c_{\alpha,\beta}(m)$ is the cocycle whose class is $\tau(m)$.  Then,  if $\F$ is an $\O$-module,  an $\A_{\tau}$-module structure on $\F$ is
given by a collection of commuting endomorphisms $\phi_{\alpha}(m)\in H^0(U_{\alpha},\End(\F))$, such that
$\phi_{\alpha}(m)-\phi_{\beta}(m)=$ multiplication by $c_{\alpha,\beta}(m) $.   The corresponding $\B_{\tau}$-module structure on $\ft(\F)$ is described
locally by 
\begin{equation}\label{eq:M local  formula}
\nabla_{m}=\nabla^{\alpha}_m+\phi_{\alpha}(m)\ .\end{equation}

\subsection{Action of $\ker(\tau)$ } Let $\Cal E'_{\tau} $ denote the restriction of \eqref{Mextension} to the kernel of $\tau$. 
\begin{equation}
\begin{CD}0@>>>\O@>{t}>>\Cal E'_{\tau}@>>>\O\otimes_{\BBB C} \ker{\tau}@>>>0\\
&&\Vert&&@VVV@VVV\\
0@>>>\O@>{t}>>\Cal E_{\tau}@>>>\O\otimes_{\BBB C} M@>>>0\\
\end{CD}
\end{equation}
The extension $\Cal E'_{\tau}$ splits canonically.  One therefore has a map of vector spaces $\ker{\tau}\to H^0(X,\Cal E'_{\tau})$,  which, in turn, gives
a map
$$\ker{\tau}\to H^0(X,\A_{\tau})\ .$$
Given an $\A_{\tau}$-module, $\F$,  one may restrict  $\nabla$  to   $\ker\tau$,  giving a map
 $\ker \tau
\morph{\nabla}\End_{\C}(\ft(\Cal F))$.  On the other hand, since
$\Cal A_{\tau}$ is commutative, there is a map $H^0(X,\A_{\tau})\to\End_{\A_{\tau}}(\F)$, and by functoriality,  a map $\End_{\A_{\tau}}(\F)\to
\End_{\B_{\tau}}(\ft(\F)) $.   The following lemma is evident from the \v Cech description \eqref{eq:M local formula}.
\begin{Lem}\label{lem:commutes}The following diagram commutes:
\begin{equation}
\begin{CD}\ker{\tau}@>>>\End_{\C}(\ft(\Cal
F))\\
@VVV@VVV\\
H^0(X,\A_{\tau})@>>>\End_{\B_{\tau}}(\ft(\F))\\
\end{CD}\end{equation}
\end{Lem}

\subsection{The algebras $\A_{\Cal C}$ and $\B_{\Cal C} $ associated to an $\O$-algebra $\Cal C$.}\label{subsect:C}  

Let $\Cal C$ be a  \qc\ sheaf of
commutative 
$\O_X$-algebras, such that
the inclusion  $\O\to\Cal C$ is injective.
Since $H^0(X,\Cal I_S)=0$,
there is
an exact sequence
\begin{equation}\label{bigkernel}
0\to H^0(X,\Cal C)\to H^0(X,\Cal C/\Cal I_S)\morph{\tau}\lie
\ .\end{equation}
The construction in subsection \ref{subsect:generalize} associates to $\tau$  
sheaves   $\A_{\tau}$ and $\B_{\tau}$,  which we denote respectively by
$\A_{\Cal C} $ and  $\B_{\Cal C}$.

\begin{Prop}\label{canonicalhomomorphism}  There is a canonical homomorphism
\begin{equation}
\A_{\Cal C}\to \Cal C
\ .\end{equation}
\end{Prop}

\begin{pf} Set $M=H^0(X,\Cal C/\Cal I_S)$. The complex $0\to\Cal C\to\Cal C/\O\to 0$ is quasi-isomorphic to $\O$.   Therefore, the  natural map  
$\O\otimes_{\BBB C} M\morph{\alpha} \Cal C/\O $ defines an extension $\Cal E\in\Ext^1(\O\otimes_{\BBB C} M,\O)
$, together with a commutative diagram
\begin{equation}
\begin{CD}0@>>>\O@>{t}>>\Cal E@>>>\O\otimes_{\BBB C} M@>>>0\label{esequence}\\
&&\Vert&&@VVV@VVV\\
0@>>>\O@>>>\Cal C@>>>\Cal C/\O@>>>0
\end{CD}
\end{equation}
It is clear that  the extension $\Cal E$ is defined by the 
natural map $M\to\lie$, so one has a canonical 
isomorphism
\begin{equation}
\Sym_{\O}(\Cal E)/\langle \
t-1 \rangle\simeq \A_{\Cal C}\ .
\end{equation}
Therefore, the map $\Sym_{\O}(\Cal E)\to\Cal C$ induced by \eqref{esequence} defines the desired homomorphism.

\end{pf}

\begin{Cor}\label{induced} $\Phi$ defines a functor from $\Cal C$ modules to $\B_C$ modules.

\end{Cor}
%
\begin{Prop}\label{sameguy}

Let $\G$ be a $\Cal C$-module, and let $m\in H^0(X,\Cal C)$.   Regarding $m$ as an element of $H^0(X,\Cal C/\Cal I_S)$,  let 
$\nabla_m$ denote the correspondng $\C$-linear endomorphism of $\ft(\Cal
G)$.  Let $m$ also denote the  $\B_C$-linear endomorphism of $\ft(\Cal
G)$ obtained by regarding $m$ as an $\A_C$-linear endomorphism of $\Cal
G$ and applying functoriality.   Then 
 $\nabla_m=m$.\end{Prop}
\begin{pf}
Since $H^0(X,\Cal 
C)$ is the kernel of the map $\tau:H^0(X,\Cal C/\Cal I_S)\to\lie$, lemma \ref{lem:commutes} implies that $\nabla_m$ acts by the image of $m$ in $\A_{\Cal
C}$. The reader can easily check that the composition
$$H^0(X,\Cal 
C)=\ker{\tau}\to H^0(X,\Cal 
A_{\Cal C})\to H^0(X,\Cal 
C)\  $$
is the identity map.  This proves the claim.

\end{pf}

\subsection{The algebras $\A_{Z}$ and $\B_{Z} $ associated to a hypersurface $Z$.} Let $Z\subset X$ be a reduced, ample Cartier
hypersurface. Let
$$
\O(*Z)=\lim_{j\to\infty}\O_X(jZ)\ .$$ 
Taking $\Cal C=\O(*Z)$ in subsection \ref{subsect:C} furnishes us with sheaves of algebras
$\A_{\O(*Z)} $ and $\B_{\O(*Z)} $. These we denote, respectively, by  
$\A_{Z} $ and    $\B_{Z} $.

For $\F$ a sheaf of $\O_X$-modules,   define
$$
\F(*Z)=\F\otimes_{\O}\O_X(*Z)\ .$$
Our 
main object of study is the  functor
$\Cal F\mapsto  \ft(\F(*Z))$ from $\O_X$-modules to $\B_{Z}$-modules.   
Denote the  $\B_{Z}$-module structure by $\nabla$.

\subsection{Local and singular parts of $\nabla$} From now on, $\F$ denotes a torsion-free
 $\O_X$-module.

\begin{Lem}\label{lem:imbedding}   Let $U$ be 
an affine open subset of $X$,  such that  $S\subset U$ and $U$ 
meets every component of $Z$.  
Then the 
natural map
\begin{equation} \ft(\F(*Z))\to 
\ft(\F(*Z)|_U)\end{equation}   is injective.
\end{Lem}
\begin{pf}  The natural map $\F(*Z)\to \F(*Z)|_U$ is 
injective,  since $\F$ is torsion-free.
  Since $\ft$ is left-exact,  the lemma is
proved.
\end{pf}

 Let $U_1,...,U_a$ be an affine open cover of $X$.    Select a \v Cech description,
$(\nabla^1_{\partial_m}+\phi_1(m),...,\nabla^a_{\partial_m}+\phi_a(m))$  of the
$\B_Z$-module structure on $\ft(\O(*Z))$.   Recall that $\nabla^i$ is a relative connection on $\P|_{U_{\alpha}\times Y} $, restricting to the standard
connection on $\P|_{(U_{\alpha}\cap S)\times Y} $.  Assume that
$U_1$ satisfies the hypotheses of lemma \ref{lem:imbedding}.  Then the lemma  implies that the  $\B_Z$-module structure on $\ft(\F(*Z))$ is
determined by
$\nabla^1_{\partial_m}+\phi_1(m)$.


\begin{Def} Let $m\in H^0(X,\O(*Z)/\Cal I_S)$ and let $U\subset X$ be an open set.   By  a lift of 
$m$ on $U$ we mean an 
element $\tilde m\in H^0(X,\O(*Z)|_U)$  such that $m$ and $\tilde m$ 
have the
same image in $H^0(X,(\O(*Z)/\Cal I_S)|_U)$.
\end{Def}

It is clear that for all $m\in H^0(X,\O(*Z)/\Cal I_S)$,  $\phi_1(m)$ is a lift of 
$m$ on $U_1$.   We therefore have the following ``local formula" for the $\B_Z$-module structure on $\ft(\F(*Z))$.

\begin{Prop}\label{decomposition} Let $U$ be as in lemma \ref{lem:imbedding}.   Then there exists a 
relative connection on $\P|_{U\times Y} $, restricting to the standard
connection on $\P|_{S\times Y}\simeq \O|_{S\times Y}$,  such that the $\B_Z$-module structure on $\ft(\F(*Z))$ is given by the following local formula:
%
%
%
%
\begin{equation}\label{localform}
\nabla_m=\nabla_{\partial_m}^U+\tilde m\ ,
\end{equation}
where $\tilde m$ is a lift of $m$ on $U$.\end{Prop}
%
%
%
%
%

For a fixed choice of $U$ and $\nabla^U$,  we refer to $\nabla_{\partial_m}^U$ and $ \tilde m$ as the local and singular parts of $\nabla_m$,
respectively.

\section{Commuting matrix differential operators}\label{sec:cmdo}

 \subsection{The algebras $\A_W$ and $\B_W$}

Let
\begin{align}\Cal N_{Z,S}&=\O(Z)/\Cal I_S\\
\Cal N_{Z}&=\O(Z)/\O
\end{align}
Then $\Cal N_Z$ is the normal sheaf of $Z$,
  which   we will regard either as a sheaf on $X$ or
$Z$ depending on the circumstance.

Note that one has an exact sequence
\begin{equation}
0\to \O_S\to  \Cal N_{Z,S}\to\Cal N_{Z}\to 0\ .
\end{equation}

  Let
$n=dim(X)$.    Let
$W\subset H^0(X,\Cal N_{Z,S})$ be an $n$-dimensional space of 
sections with the following properties:\medskip
\begin{align}\label{nicew}
&\kern 70pt \text{  $W$ is isomorphic to its image in $H^0(X,\Cal 
N_{Z})=H^0(Z,\Cal 
N_{Z})$.}\\
&\kern 70pt\text{  The image of $W$ in
$H^0(Z,\Cal N_{Z})$ is   basepoint-free.}  \end{align}
Then $W$ defines a finite morphism
\begin{equation}
\map:Z\to\BBB P^{n-1}=\BBB P(W^*)\ .
\end{equation}

 There is a  natural map $W\to\Frak g$, factoring through the 
inclusion $W\subset H^0(X,\O(*Z)/\Cal I_S)$.

\begin{Def} Let $\A_W$  and $\B_W$ denote the sheaves of algebras 
associated to the natural map $W\to\Frak g$.  Note that $\A_W\subset\A_Z$ 
and $\B_W\subset\B_Z$.\end{Def}

\subsection{Filtration} From now on we assume $\F$ is rank-one,  torsion free,  and  locally 
free in a neighborhood of $Z$. Denote $\F(k Z)$ by $\F(k)$.  Then $\{\F(k)\} 
$ is  a filtration of $ \F(*Z)$.   The associated graded sheaf,
$\oplus_k \F(k)/\F(k-1) $ is  a sheaf of 
$\O_X\otimes_{\C}\Sym(W)$-modules.   Now $\{\ft(\F(k))\} $ is  a 
filtration of $ \ft(\F(*Z))$ and we
have an imbedding
\begin{equation}gr(\ft(\F(*Z)))\subset \ft(gr(\F(*Z)))\ .\end{equation}
There is also a filtration on $\B_{W} $,  with respect to which
\begin{equation}gr(\B_W)\simeq\O_{\gp}\otimes_{C}\Sym(W)\ .\end{equation}

\begin{Lem} The action of $\B_W$ on $\ft(\F(*Z))$   respects the 
filtration.   Moreover,   the action of $gr(\B_W)$ on 
$gr(\ft(\F(*Z)))$
is compatible with the action of $\Sym(W)$ on $\ft(gr(\F(*Z)))$ 
induced by the action of $\Sym(W)$ on $gr(\F(*Z))$.
\end{Lem}

\begin{pf}    Let $m\in W$.  By proposition \ref{decomposition}, the 
singular part of $\nabla_m$ is given,  locally in
$X$,   by multiplication by a rational  function with poles on $Z$, 
whose  class  in $\Cal N_{Z,S} $ is  $m$ itself.   This proves the 
lemma.
\end{pf}

\subsection{Local freeness}
Next we give a criterion for $\ft(\F(*Z))$  to be a locally free 
$\B_W$-module on a dense open subset of $\gp$.

\begin{Thm}\label{Locally Free at the generic point} Assume
that for a
general line bundle $\Cal M$ on $X$ of degree $0$,
\begin{equation}H^k(X,\Cal M\otimes \F(j))=0\end{equation}for  all 
$j$  and $0< k<n$.

Then $\ft(\F(*Z))$ is locally free as a $\B_W$-module
at a general  point of
$\gp$.   Moreover,
\begin{equation}\label{rank}
rank_{\B_W}(\ft(\F(*Z)))=deg(\Cal N_Z)=deg(\map)\ .
\end{equation}

\end{Thm}

\begin{pf}  It suffices to prove that for general $\Cal M$,  the 
fiber of $gr (\ft(\F(*Z)))$
at $\Cal M$ is free as a $\Sym(W)$-module.

  First consider the case $n\ge 2$.
By the hypothesis,  with $k=1$, the fiber is
\begin{equation}\label{assocgraded}
\Oplus_{j} H^0(Z,\Cal M\otimes \F(j)|_Z)\ .
\end{equation}
Let
\begin{equation}
\Cal S_{\Cal M}=\map_*(\Cal M\otimes \F|_Z)\ .
\end{equation}
Then $\Cal S_{\Cal M}$ is a vector bundle on $\BBB P^{n-1}$.  Moreover,
  \begin{equation}
H^k(\BBB P^{n-1},\Cal S_{\Cal M}(j))=0
\end{equation}
for all $j$  and $0< k< n-1$.
By Horrock's criterion \cite{OSS},  $\Cal S_{\Cal M}$ decomposes as 
sum of line bundles.     Then
$\Oplus_j H^0(\BBB P^{n-1},\Cal S_{\Cal
M}(j))$ is free as a $\Sym(W)$-module. (Note that the sum is over all 
integers $j$,  but is bounded below.)  This proves the result when 
$n>2$.

When $n=1$,   the argument is slightly different.   $Z$ is now a sum 
of smooth points,  $p_1,...,p_d$,   and
$\Oplus_{j} H^0(Z,\Cal M\otimes \F(j)|_Z)=\oplus_j \C[\xi,\xi^{-1}]$,
where $\xi$ is a basis for $W$.   On the one hand,   this is not a 
free $ \C[\xi]$-module,  but on the  other hand,  the fiber of
$gr (\Phi(\F(*Z)))$ is a proper submodule,   isomorphic to $\oplus_k \C[\xi]$.

\end{pf}

As a corollary, we obtain the first theorem announced in the introduction.

\begin{Cor}\label{cor:CM}   Assume $X$ is Cohen-Macauley and $\F$ is locally free. 
Let  $\cal L$  be an ample line bundle on $X$.   Let $Z'\subset X$ be 
an ample
Cartier hypersurface. Then there exist integers
$a$ and
$b$ such that if $Z$ is a  reduced hypersurface algebraically 
equivalent to $bZ'$ and $W\subset H^0(X,\Cal N_{Z,S})$ is an 
$n$-dimensional space such
that properties
\eqref{nicew} hold,   then $\ft(\Cal L^a\otimes\F(*Z)) $ is 
locally free as a $\B_W$-module
at a general  point of
$\gp$.
\end{Cor}

\begin{pf}   First choose $a$ so that   for all $k>0$,  and general 
degree-zero line bundle $\Cal M $, $H^k(X, \Cal L^a\otimes\Cal
M\otimes\F)=0$.   Then choose positive integer $b$ such that $H^k(X, 
\Cal L^a\otimes\Cal M\otimes\F(jZ'))=0$ for $
|j|\ge b$ and $0<k<n $.   The hypothesis of theorem \ref{Locally Free at the
generic point} is now satisfied,  with $\F$ replaced by $\Cal 
L^a\otimes \F $.\end{pf}

\begin{Rem} The hypothesis of theorem \ref{Locally Free at the 
generic point} makes no
reference to the subscheme $S$.   Thus one may regard $S$ as a 
variable,  and fix a choice of $S$ at the end.
\end{Rem}

\subsection{Imbedding of rings } 
As noted in the introduction, upon choosing a 
 $\B_W$-basis for $\ft(\F(*Z))$ over an open set
$\Cal U\subset \gp$,  the functorially defined imbedding
\begin{equation}\label{intrinsicimbedding} H^0(X,\O(*Z))\to 
\End_{\B_{Z}}(\ft(\F(*Z)))\ \end{equation}
 becomes an 
imbedding 
\begin{align}\label{imbedding}H^0(X,\O(*Z))&\to H^0(\Cal U,M_d(\B_W))\\
f&\mapsto L_f\ .\notag\end{align}
into the ring of $d\times d$ matrices with entries in $H^0(\Cal U, \B_W)$. 

\subsection{Filtered basis} 
Assume henceforth that $\F$ satisfies the conditions 
of theorem \ref{Locally Free at the
generic point}.

\begin{Def} Let $d$ denote the rank of
$\ft(\F(*Z))$ as a $\B_W$-module. Let 
$(c_1,...,c_d)$ be the degrees of a set of homogeneous generators 
of the fiber of $gr
(\ft(\F(*Z)))$ as a $\Sym(W)$-module, listed in nondecreasing 
order. By filtered basis of $\ft(\F(*Z))$ over an open set
$\Cal U\subset \gp$,  we mean a column vector 
$(\psi_1,...,\psi_d)^T$  such that for all $i$, 
$\psi_i\in\Gamma(\Cal U, \ft(\F(c_i))$, and the classes
$$
\psi_i\  mod\ \ft(\F(c_i-1))$$ are a basis for
$gr(\ft(\F(*Z)))|_{\Cal U}$.\end{Def}

Unless otherwise stated,  we agree henceforth to work only with filtered bases.
We then introduce a filtration on  the matrices over $\B_W$, in such a way that  
\eqref{imbedding} 
is a
homomorphism of filtered algebras.
\begin{Def}\label{matrixfiltration} 
Given a positive integer $j$,  let $\mat j$ denote  the  sheaf of
 $j\times j$ matrices with entries in 
$\B_W$. 
Given a 
nondecreasing sequence of integers $\vec r=(r_1,...,r_j)$, define
the $\vec r$-filtration on $\mat j$  by
\begin{equation}\label{rfiltration}
{\mat j}_{\vec r,k}=
\{ \  (L_{a,b})\ |\ ord(L_{a,b})\le k+r_a-r_b\ \}\ .
\end{equation}
\end{Def}

\begin{Prop}\label{filteredimbedding}  If ${\mat d}$ is endowed with
the $(c_1,...,c_d)$-filtration and if  the homomorphism
\eqref{imbedding} is defined with respect to a filtered basis, then 
it is a   homomorphism of filtered algebras.  Moreover,  the induced 
map of
associated graded algebras is an imbedding.\end{Prop}

\begin{Rem}

Let $\Cal S_{\Cal M}$ be the vector
bundle on
$\BBB P^{n-1}$ introduced in the proof of theorem
\ref{Locally Free at the generic point}.  If $n>1$ then 
$-c_1,...,-c_d$ are the Chern classes appearing in the
decomposition of $\Cal S_{\Cal M}$ as a sum of line bundles.  \end{Rem}

\subsection{Zero Curvature Equations;  Flows, I} 
There is a map 
\begin{align}H^0(X,\O(*Z)/\Cal I_S)&\to H^0(\Cal U,\mat d)\notag\\
m&\mapsto L_m\end{align}
defined by letting $\nabla_m$ act on the basis:
\begin{equation}
  \nabla_m(\psi)=L_m\psi\ .
\end{equation}
For $f\in H^0(X,\O(*Z))$, there are then two definitions of $L_f$.  But by proposition \ref{sameguy}, the two definitions coincide.
%
\begin{Prop}\label{prop:zero curvature}

For all $m,n\in H^0(X,\O(*Z)/\Cal I_S)$,  
one has  the {\em zero-curvature equation}
\begin{equation}
  \label{eq:zerocurvature}
  \partial_m(L_n)-\partial_n(L_m)=[L_m,L_n]\ ,
\end{equation}
where vector fields are acting on differential operators by acting on the 
coefficients.  In particular,   for all $f\in H^0(X,\O(*Z))$,
\begin{equation}\label{flow1}
\partial_m(L_f)=[L_m,L_f]\ .
\end{equation}
\end{Prop}
\begin{pf} The first assertion follows immediately from the fact that $\nabla_m$ commutes with $\nabla_n$.   For the second assertion,  set
$n=f$ and note that $H^0(X,\O(*Z))$ is the kernel of the map $H^0(X,\O(*Z)/\Cal I_S)\to H^1(X,\Cal I_S)$.\end{pf}

Our problem now is to choose the basis in such a way that the  matrices $\{L_m\ |\ m\in H^0(X,\O(*Z)/\Cal I_S)\}$ are given explicitly in terms
of the algebra of matrices $\{L_f\ |\ f\in H^0(X,\O(*Z)\}$.   The first step is to choose a good projective basis, that is to say a 
choice of $\psi_1,...,\psi_d$  up to multiplicative factors in $\O_Y^*$.    Then we find many ways to lift  a
good projective basis to a good basis,   leading to  equation
\eqref{eq:flows zero} and its various modifications.

\section{Good Basis,  Projectively}\label{sec:good  basis}

\subsection{Laurent expansion}

Let $Q$ be a smooth point of $X$.  Let  $\hat{\O}_Q $  
denote the formal completion of the stalk of
$\O_X$ at $Q$, regarded as a sheaf of $\O_X$-modules.
Let $\Frak m_Q $ denote the maximal ideal in $\hat{\O}_Q $.
\begin{Def}  For all integers $a$ and $k$,  let\begin{align}
\ft(\F(k)\otimes\Frak m_Q^a)&=\lim_{\overset \leftarrow \ell}
\ft(\F(k)\otimes\Frak m_Q^a/\Frak m_i^{\ell})\ .\\
\ft(\F(*Z)\otimes\Frak m_Q^a)&=\lim_{\overset \rightarrow k}
\ft(\F(k)\otimes\Frak m_Q^a)\ .
\end{align}
\end{Def}

\begin{Prop}   For every smooth point $Q\in Z$, 
  there is a natural imbedding of $\B_{Z}$-modules
  \begin{equation}
    \label{eq:laurentexpansion}
    \ft(\F(*Z))\to\ft(\F(*Z)\otimes\hat{\O}_Q) \ .
  \end{equation}
\end{Prop}

\subsection{Gaussian elimination}
Fix a point $P\in \BBB P^{n-1}$, such that $\map$ is \'etale in a 
neighborhood of $\map^{-1}(P)$.  Fix a basis
$\xi,\eta_1,...,\eta_{n-1}$ for $W$ such that
\begin{align}\xi(P)&\ne 0\\
\eta_i(P)&=0\ ,\ i=1,...,{n-1}\ .
\end{align}

 Let $U$ and $\nabla^U$ be as in  
proposition \ref{decomposition},  and assume also that 
$U$   contains $ \map^{-1}(P)$. 
One obtains 
 rational  functions $z,w_1,...,w_{n-1}$ on $X$ by taking the singular part of $\nabla$,  i.e.,  
\begin{align}\nabla_{\xi}&=\nabla_{\partial_{\xi}}^{ U}+\frac 1z\label{coordinates1}\\
\nabla_{{\eta_i}}&=\nabla_{{\partial_{\eta_i}}}^{U}+\frac{w_i}z\label{coordinates2}
\end{align}
Then at every point in $\map^{-1}(P)$, $(z,w_1,...,w_{n-1})$ 
generate the maximal 
ideal   and  $z$ generates the ideal of  $Z$.

Now let $\vec Q=(Q_1,...,Q_j)$ be a list of distinct points in
 $\map^{-1}(P)$ and let $\vec r=(r_1,...,r_j)$ be a 
nondecreasing sequence of
integers.    Define
\begin{equation}
 \zeta_{\vec r}=\pmatrix\begin{array}{llll}\frac{1}{z^{r_1}} &0  &
\dots&0\\0&\frac{1}{z^{r_2}}&
\dots&0\\
\vdots&
\vdots&
\ddots&
\vdots\\0&0&
\dots&\frac{1}{z^{r_j}}
\end{array}\endpmatrix\ . 
\end{equation}

\begin{Def} For any integer $\ell$, let 
$\Frak m_{\vec Q,\vec r}^{(\ell)} $
 denote the sheaf of $j\times j$ 
matrices whose entry in 
row $a$ 
column $b$ is the ideal $\Frak m_{Q_a}^{r_a-r_b+\ell-\delta_{a,b}}$,
where $\Frak m_Q^r=\hat {\O}_Q $ for $r\le 0$.
\end{Def}

For every point $Q$ in $X$,  the structure sheaf of $Q$ transforms 
to a line 
bundle on $\gp$.   Denote this line bundle by $\Cal L_Q$.    Fix a
trivialization of $\Cal F$ in a neighborhood of $\map^{-1}(P)$. 
Then for every point $Q\in\map^{-1}(P)$ we
 have a
morphism
\begin{equation}\label{modmax}
\ft(\F\otimes\hat\O_{Q})\to\Cal L_{Q} 
\end{equation} 
with kernel $\ft(\F\otimes\Frak m_{Q}) $.    Refer to \eqref{modmax} as
``evaluation at $Q$.''
Suppose then that one has a $j\times j$ matrix 
\begin{equation}\label{chi}
  \chi=(\chi_{a,b})\end{equation}
such that for all $a,b$,
\begin{equation}\chi_{a,b}\in\ft(\F(r_a)\otimes\hat\O_{Q_b})\ .
\end{equation}
 Set
\begin{equation}
  \label{eq:tildechi}
  \tilde\chi=\zeta_{\vec r}^{-1}\chi \ .
\end{equation}
Then the $(a,b)^{th}$ 
entry  of $\tilde\chi$ is a section of $\ft(\F\otimes\hat \O_{Q_b})$,  so it 
makes sense   to evaluate the entries of $\tilde\chi$.  Thus we have a matrix
\begin{equation}
  \label{eq:evaluate chi}
  \tilde\chi|_{\vec Q}
\end{equation}
such that each row is a section of the vector bundle $\Oplus_b \Cal L_{Q_b}$.
It then makes sense to ask whether $\tilde\chi|_{\vec Q}$ is invertible.

\begin{Prop}   If
$\tilde\chi|_{\vec Q}$ is invertible then there exists a
permutation matrix $\sigma$ and an invertible degree-$0$ matrix
$L\in\mat j$, such that $\widetilde{L\chi\sigma}$ is a section
of $\ft(\F\otimes\Frak m_{\vec Q,\vec r}^{(1)})$.
\end{Prop}
\begin{pf}
It is clear that any composition of moves of the following two types
will transform $\chi$ to a matrix of the form $L\chi\sigma$:

\medskip
1.  reordering the columns,

\medskip
2. replacing $\chi_{a_2}$ by $\chi_{a_2}+L_0 \chi_{a_1} $ if $L_0$ is a 
section of $\B_W$ such that $$ord(L_0)\le r_{a_2}-r_{a_1} .$$
\medskip
Note that a move of type 2 replaces $\tilde\chi_{a_2}$ by 
$\tilde\chi_{a_2}+z^{r_{a_2}-r_{a_1}}L_0 \tilde\chi_{a_1} $.   It is 
then clear from
equations \eqref{coordinates1} and \eqref{coordinates2} that if 
$\tilde\chi_{a_1,a_1}|_{\vec Q}$ is invertible,  then for all $a_2\ne a_1$, 
a move of type 2 can be 
used to move
$\tilde\chi_{a_2,a_1}$ into $\ft(\F\otimes\Frak 
m_{Q_{a_1}}^{r_{a_2}-r_{a_1}+1} )$.

\end{pf}
\begin{Rem}
  If $\tilde\chi$ is a section of 
$\ft(\F\otimes\Frak m_{\vec Q,\vec r}^{(1)})$, then $\tilde\chi|_{\vec Q}$ 
is upper triangular,  and is diagonal within each block where $r_a=r_b$.
\end{Rem}

From now on we  let $\vec Q=(Q_1,...,Q_d)$ denote an ordering of $\map^{-1}(P)$.   
 Let $\psi=(\psi_1,...,\psi_d)^T$ be a filtered
basis.  Let $\hat\psi$ denote the $d\times d$ matrix obtained by 
expanding each row of $\psi$ at each point $Q_b$.
Then the $(a,b)^{th}$ entry of $\hat\psi$ is a section of 
$\ft(\F(c_a)\otimes\hat\O_{Q_b})$.  One then has the matrix
\begin{equation}
\tilde\psi=\zeta_{\vec c}^{-1}\hat{\psi}\ .
\end{equation}
The key point of this section is the following proposition.
\begin{Prop}
The matrix $\tilde\psi|_{\vec Q}$ is invertible.
\end{Prop}
  \begin{pf}
     Consider the vector bundle on $\BBB P^{n-1}\times \gp$,
\begin{equation}\Cal S=(\map\times id)_*(\pi_Z^*(\F|_Z)\otimes\P)\ 
.\end{equation}
Let $gr(\psi)=(gr(\psi)_1,...,\gr(\psi_d))^T$ be the graded basis 
associated to $\psi$.  Then $gr(\psi)$ establishes an isomorphism
\begin{equation}\label{sumoflines}\Cal S|_{\BBB P^{n-1}\times\Cal 
U}\simeq \Oplus_{i=1}^d\O(-c_i)\ .\end{equation}
That is,  for all $i$,
\begin{equation}gr(\psi_i)\in\Gamma (\BBB P^{n-1}\times\Cal U,\Cal 
S(c_i))\ ,\end{equation}
and with respect to the isomorphism  \eqref{sumoflines}, $gr(\psi)$ 
is the identity matrix.
Let
\begin{equation}\label{xi}\Xi=\pmatrix\begin{array}{llll}{\xi^{c_1}} & &
&\\&{\xi^{c_2}}&
&\\
&
&
\ddots&
\\&&
&{\xi^{c_d}}
\end{array}\endpmatrix\ .
\end{equation}
Now $\xi^{-c_i}$ trivializes $\O(-c_i)$ at $P$. Then $gr(\psi)_i$ may 
be evaluated along $\{P\}\times\Cal U$ to produce a section
\begin{equation}
\xi^{-c_i}gr(\psi_i)|_P\in\Gamma (\{P\}\times\Cal U,\Cal 
S|_{\{P\}\times \gp})\ .\end{equation}
On the other hand,
\begin{equation}\label{othersumoflines}\Cal S|_{\{P\}\times 
\gp}\simeq \Oplus_{j=1}^d\Cal L_j\ .\end{equation}
It is easy to see that with respect to the isomorphism \eqref{othersumoflines},
\begin{equation}
\tilde\psi|_{\vec Q}=\Xi^{-1}\  gr(\psi)|_P\ .
\end{equation}
This proves the proposition.
  \end{pf}

  \begin{Def}[Good Basis]   
    Say that a filtered basis $\psi$ is {\em good} with respect 
to $\vec Q$ if $\tilde\psi$ is a section
of $\ft(\F\otimes\Frak m_{\vec Q,\vec c}^{(1)})$.   In other words,   
$\psi$ is good if, after trivializing $\F$ at each point
of $\map^{-1}(P)$, its matrix of Laurent expansions at $\vec Q$ 
has the form
\begin{equation}\kern-40pt
\pmatrix\begin{array}{llll}\frac{1}{z^{c_1}} &0  &
\dots&0\\0&\frac{1}{z^{c_2}}&
\dots&0\\
\vdots&
\vdots&
\ddots&
\vdots\\0&0&
\dots&\frac{1}{z^{c_d}}
\end{array}\endpmatrix
\bigg(
\begin{array}{lll}\text{invertible}\\\text{diagonal}\\ \text{matrix}\end{array}
+
\pmatrix\begin{array}{lllll}\Frak m_{Q_1} &\Frak m_{Q_1}^{c_1-c_2+1} &
\dots&\Frak m_{Q_1}^{c_1-c_d+1}\\\Frak m_{Q_2}^{c_2-c_1+1}&\Frak m_{Q_2}&
\dots&\Frak m_{Q_2}^{c_2-c_d+1}\\
\vdots&
\vdots&
\ddots&
\vdots\\\Frak m_{Q_d}^{c_d-c_1+1}&\Frak m_{Q_d}^{c_d-c_2+1}&
\dots&\Frak m_{Q_d}
\end{array}\endpmatrix\bigg)\kern30pt
\end{equation}
  \end{Def}
  \begin{Cor}
    \label{goodbasis}  Given any filtered basis $\psi$,   there exists an 
invertible  degree-$0$ matrix
$L\in\mat d$ such that $L\psi$ is good with respect to some ordering of 
$\map^{-1}(P)$.
  \end{Cor}

\section{Microdifferential operators}\label{sec:microdiff}

\subsection{Microdifferential operators}Imbed $\B_W$ into a larger ring $\pB$ in the following standard way. 
Sections of $\pB$ are formal expressions
\begin{equation}\label{mdop}L=\sum_{i=-\infty}^N L_i\xi^{i}\ ,\end{equation}
where $L_i=L_i(\eta_1,...,\eta_{n-1})$ are sections of the subsheaf 
of $\B_W$ generated by $\eta_1,...,\eta_{n-1}$,  and for 
$f\in\O_{\gp}$,
\begin{equation}\xi^{-1}f=\sum_{i=-\infty}^{-1} 
(-1)^{i+1}\partial_{\xi}^{-i-1}(f)\xi^{i}\ .\end{equation}
Note that if $\eta$ is any $\C$-linear combination of 
$\eta_1,...,\eta_{n-1}$  and $a$ is a nonzero constant, then
\begin{equation}{(a\xi+\eta)}^{-1}=\sum_{i=-\infty}^{-1} 
a^{i}{\eta}^{-i-1}\xi^{i}\ \end{equation}
in $\pB$, so that indeed the definition depends only on $P$.

Say 
that the section $L$ in \eqref{mdop} has order at most $k$ if 
$ord(L_i)+i\le k $
for all
$i$.   Let \dsp{\pB_k} denote the subsheaf consisting 
 elements of order at most 
$k$ and let \dsp{\pBf=\lim_{\overset k\to}\pB_k}.
\begin{Prop}   The action of $\B_W$ on $\ft(\F(*Z))$ extends 
canonically to a filtered action of $\pBf$ on $\Oplus_j
\ft(\F(*Z)\otimes\hat\O_j)$.   If  $\psi_1,...,\psi_d$ is a 
filtered basis for  $\ft(\F(*Z))$ as a $\B_W$-module,  then
$\hat\psi_1,...,\hat\psi_d$ is a filtered basis for  $\Oplus_j
\ft(\F(*Z)\otimes\hat\O_j)$ as a $\pBf$-module.
\end{Prop}

The proof is easy and well-known for curves, and works equally well 
in any dimension.

Given $\sum L_i\xi^{i} 
\in\pB$, define its
differential-operator part and negative part, respectively,   by
\begin{equation}(\sum L_i\xi^{i})_{+} =\sum_{i\ge 
0}L_i\xi^{i}\ ,\  (L_i\xi^{i})_-= \sum_{i<0} L_i\xi^{i} \ 
.\end{equation}
Extend the definition of $(\cdot)_{\pm} $ to matrices,  entry by entry.

\subsection{Flows, II}
Fix a  basis $\psi$. The  local and singular parts of $\nabla_m$
separately act on $\Oplus_j
\ft(\F(*Z)\otimes\hat\O_j)$.   Thus one obtains matrices 
$L_m^U$ and $L_{\tilde m} $  over $\pBf$, defined by
\begin{equation}
\nabla_{\partial_m}^U\hat\psi=L_m^U\hat\psi\ ,
\end{equation}
\begin{equation}
{\tilde m}\hat\psi=L_{\tilde m}\hat\psi\ .
\end{equation}

With respect to the $\vec c$-filtration,  $L_m^U$ has order $0$ while  
$L_{\tilde m} $, has the same order as $m$.  Clearly,   $L_m= L_m^U+{L_{\tilde 
m}}$.   Thus
\begin{align}\label{psdecomp}L_m&=(L_m^U)_{+}+({L_{\tilde 
m}})_{+}\\\label{negpart}(L_m^U)_{-}&=-({L_{\tilde 
m}})_{-}\ .\end{align}

One now has,  for all $h\in H^0(X,\O(*Z))$,
\begin{equation}\label{flow2}
\partial_m(L_h)=[(L_m^U)_++(L_{\tilde m})_+,L_h]\ .
\end{equation}

Notice that the lift
 $\tilde m$ is a rational function on $X$.  Thus $(L_{\tilde m})_+$ is the term appearing in equation \eqref{eq:flows zero}.  Let us call $(L_m^U)_+$ the {\em
local term}.   The local term is the  
basis-dependent part.
If the basis is good,  then  the local term,   being of order $0$,   is 
block lower-triangular,  with $0^{th}$ order operators,  (i.e. 
sections of $\O$) in the
diagonal  blocks.   Fortunately,   we can say more.

\begin{Lem}\label{connection1forms} Assume the basis $\psi$ is 
good.  Then $(L_m^U)_{+}$ is
a diagonal matrix with entries in $\O_Y$.   In fact, if we let
$A$ denote the diagonal part of $\tilde\psi$, then
\begin{equation}\label{logdiff}(L_m^U)_{+}=
\nabla_{\partial_m}^U(A|_{\vec Q})A|_{\vec Q}^{-1}\ .
\end{equation}
\end{Lem}
\begin{pf}  We make use of the following facts:   
Let $r\in\pBf$ be 
an element of order $\le k$.   Let $a$ be an integer and let
$\chi\in\ft(\F(a)\otimes\hat\O_j)$  be a section such that 
$z^a\chi(Q_j)$ is a nonvanishing section of $\Cal L_j$. Then
\medskip

A.\ \ \ $\forall\mu\forall\nu\forall c$\ 
($r(\chi)\in\ft(\F(a+k)\otimes\Frak m_j^{\nu})\ \Rightarrow\ $

$r(\ft(\F(c)\otimes\Frak 
m_j^{\mu}))\subset\ft(\F(c+k)\otimes\Frak m_j^{\nu+\mu})\ )$.
\medskip

B.\ \ \  $r(\chi)\in\ft(\F(a+k)\otimes\Frak m_j^{k+1})\ 
\Leftrightarrow\ (r)_{+}=0$.
\medskip

  To simplify notation,   set
$R=L_m^U$.    We must prove that for $i\ne j$,
$R_{i,j}(\hat\psi_{j,j})\in \ft(\F(c_i)\otimes\Frak 
m_j^{c_i-c_j+1})$.  Let $k$ and  $\ell$ be  integers such that
\medskip

1. $i\ne j\ \Rightarrow$
$R_{i,j}(\hat\psi_{j,j})\in \ft(\F(c_i)\otimes\Frak 
m_j^{\text{inf}(k,c_i-c_j+1)})$.

\medskip

2. $i\ne j$ and $ c_i- c_j>\ell\ \Rightarrow$
$R_{i,j}(\hat\psi_{j,j})\in \ft(\F(c_i)\otimes\Frak 
m_j^{\text{inf}(k+1,c_i-c_j+1)})$.

\medskip
Note that 1. holds for $k=0$ and 2. holds for  $\ell=c_d-c_1$.

Fix indices $i,j$ such that $i\ne j$ and $c_i-c_j=\ell$. We have

\begin{equation}
R_{i,j}(\hat\psi_{j,j})=\nabla_{\partial_m}^U(\hat\psi_{i,j})-\sum_{\overset{c_i-c_p\le\ell}{p\ne
j}}R_{i,p}(\hat\psi_{p,j})-\sum_{c_i-c_p>\ell}R_{i,p}(\hat\psi_{p,j})
\end{equation}
  The first and third terms belong to $\ft(\F(i)\otimes\Frak 
m_j^{\text{inf}(k+1,c_i-c_j+1)})$.
By fact A,   if
$c_i-c_p\le\ell$ and
$ p\ne j $,  $R_{i,p}(\hat\psi_{p,j})\in \ft(\F(c_i)\otimes\Frak 
m_j^{\text{inf}(k,c_i-c_p+1)+c_p-c_j+1})$.  Now $c_p-c_j=
\ell-(c_i-c_p)\ge 0$, so
$\text{inf}(k,c_i-c_p+1)+c_p-c_j+1\ge \text{inf}(k+1,c_i-c_j+1)$. 
Thus,  by descending induction on $\ell$, 2. holds for
all $\ell$.     Then 1. holds for all $k$ by
induction on $k$.    We then have $(R_{i,j})_{+}=0 $ for 
$i\ne j$,  by fact B.

To see that   formula \eqref{logdiff} gives the diagonal entries, 
recall that the diagonal entries of $(R)_{+}$ are sections 
of $\O$.   Then
  evaluate the equation
\begin{equation}
\kern -30pt 
(R_{i,i})_{+}\hat\psi_{i,i}=\nabla_{\partial_m}^U(\hat\psi_{i,i})-\sum_{p\ne
i}R_{i,p}(\hat\psi_{p,i})-(R_{i,i})_-(\hat\psi_{i,i})\end{equation}  
at $Q_i$.
\end{pf}

\section{The Flows}\label{sec:flows}

To obtain \eqref{eq:flows zero},  choose
$S=\map^{-1}(P)$.
Then the line bundles $\Cal L_j=\Cal L_{Q_j}$ are canonically 
trivial.   Thus  the entries of
$\tilde\psi|_{\vec Q}$ are sections of
$\O_{\gp}$, so the following definition makes sense.

\begin{Def}Say that a good basis 
$\psi$ is canonically normalized if the diagonal entries 
$\tilde\psi|_{\vec Q}$ are
constant.\end{Def}

\begin{Thm}\label{draconian} Assume 
$\psi$ is canonically normalized.  Then    for all $m\in
H^0(X,\O(*Z)/\Cal I_S)$ there exist $f,g\in H^0(X,\O(*Z))$ such 
that for  all $h\in H^0(X,\O(*Z))$,
\begin{equation}\label{flow3}\partial_m(L_h)=[(L_g^{-1}L_f)_{+},L_h]\ ,
\end{equation}
Moreover,   \dsp{\frac fg} is a  lift of $m$ in a 
neighborhood of $\map^{-1}(P)$.
\end{Thm}

\begin{pf} The field of rational functions on $X$ is the field of 
fractions  of   $H^0(X,\O(*Z))$.   Thus,  the local lift $\tilde m$ is of the 
form
\dsp{ f/g} for some $f,g\in H^0(X,\O(*Z))$.  In light of
equation
\eqref{flow2}, all that remains is to show that $(L_m^U)_{+}=0 $.
The relative connection $\nabla^U$ restricts on each $\Cal L_j$ to the 
standard connection on
the trivial bundle.  Thus it annihilates constants.     Then,  we are
done, by  lemma \ref{connection1forms}.
\end{pf}

\section{Flows associated to other normalizations}\label{sec:more flows} 

If one considers equations \eqref{psdecomp} and \eqref{flow2} without any 
prejudice as to the outcome,
one sees that the task of a basis is not necessarily to make the local term 
equal to zero,  but simply to make the local term a (preferably explicit) 
function of
$L_{\tilde m}$.    Here we present a procedure to do this, the idea being to rigidify  the Poincar\'e line bundle at as few points as possible.
(If a point $Q_j$ is omitted 
from $S$,  then  $\Cal L_j$ is no longer trivial, so the canonical normalization no
longer makes sense.)  The procedure is randomized,  so 
that roughly speaking,   the number of
possible flows that may be obtained by it grows  exponentially 
with the number of different degrees appearing in a 
good basis.

   Given a section of $\pBf$,   \dsp{L=\sum L_i\xi^{i}}, of order at 
most $N$,  define its $N^{th}$-order $\xi$-symbol by
\begin{equation}\label{psymbol}\sigma_{\xi}^N(L)= L_N\xi^{N}\ .\end{equation}
Note that $L_N$ belongs to  $\O_{\gp}$.   Define  the $N^{th}$-order 
$\xi$-symbol of a matrix  $L\in{\matp d} $ by
\begin{equation}\label{matsymbol}\sigma_{\xi}^N(L)_{i,j}=\sigma_{\xi}^{N+c_i-c_j}(L_{i,j})\ 
.\end{equation}
Then
\begin{equation}\xi^{-N} \Xi^{-1}   \sigma_{\xi}^N(L)\ \Xi\end{equation}
is  a matrix with entries in $\O$,   where $\Xi$ is defined by \ref{xi}.

\begin{Lem}\dsp{\Xi^{-1}\sigma_{\xi}^0(L_m^U)\ 
\Xi=\nabla_{\partial_m}^U(\tilde\psi|_{\vec Q})(\tilde\psi|_{\vec Q})^{-1}}.\end{Lem}

\begin{pf} If $a$,  $b$ and $k$ are integers such that $a+b\le k $, 
$a\ge 0$ and    $b<k$, then  $a+b<k$ or $a>0$.   Either way,  for all 
$i$, $j$ and $\ell$,
\begin{equation}\eta_i^a\xi^b\ft(\F(\ell)\otimes\hat\O_j)\subset\ft(\F(\ell+k)\otimes\Frak 
m_j)\ .\end{equation}
Therefore,
\begin{equation}\zeta^{-1}\sigma_{\xi}^0(L_m^U)(\hat\psi)|_{\vec Q}=\zeta^{-1} 
L_m^U(\hat\psi)|_{\vec Q}\ .\end{equation}
Furthermore,
\begin{align}\zeta^{-1}\sigma_{\xi}^0(L_m^U)(\hat\psi)|_{\vec Q}&=
\zeta^{-1}\sigma_{\xi}^0(L_m^U)\zeta(\tilde\psi)|_{\vec Q}\notag\\
&=\Xi^{-1}\sigma_{\xi}^0(L_m^U)\ \Xi(\tilde\psi|_{\vec Q})\ .\end{align}
On the other hand,
\begin{align}\zeta^{-1} 
L_m^U(\hat\psi)|_{\vec Q}&=\zeta^{-1}\nabla_{\partial_m}^U(\hat\psi)|_{\vec Q}\notag\\
&=\nabla_{\partial_m}^U\zeta^{-1}(\hat\psi)|_{\vec Q}\notag\\
&=\nabla_{\partial_m}^U(\tilde\psi|_{\vec Q})\ .\end{align}
Since $\Xi^{-1}\sigma_{\xi}^0(L_m^U)\ \Xi$ is a matrix over $\O$, 
the lemma is proved.
\end{pf}

The normalization problem  is now the following:

\begin{Prob}\label{nproblem3}
Denote $\nabla_{\partial_m}^U(\cdot)$ by $(\cdot)'$.   Let 
$R=(\tilde\psi|_{\vec Q})'(\tilde\psi|_{\vec Q})^{-1}$.
Choose the diagonal entries of $\tilde\psi|_{\vec Q}$ so as to 
recover the diagonal entries of $R$
from the upper triangular entries of $R$.\end{Prob}

 Before giving the 
procedure,   some illustrative examples are in order.

\begin{Ex}\label{flowexample}
\medskip

If ${\vec Q}=(Q_1,Q_2)$ and $Q_1\in S$,   then we may partially normalize $\psi$ by
\dsp{\tilde\psi|_{\vec Q}=\begin{pmatrix}1&a\\0&b\end{pmatrix}},  with $b$ 
to be determined.    Then
\dsp{R=\begin{pmatrix}0&\frac{a'}b\\0&\frac{b'}b\end{pmatrix}}. 
Recall that $a$ and $b$ are sections of $\Cal L_2$. After declaring 
$Q_2$ to be a member of
$S$,  we can proceed as in   theorem
\ref{draconian} to set
$b=1$.

However, if
$a$ is an invertible section of $\Cal L_2$,  a different choice 
presents itself,  namely   $b=a$.
How does this alter the flow \eqref{flow3}?     Let 
\begin{equation}
L=(L_g)^{-1}L_f \ . 
\end{equation}
Let $L_{i,j}$ denote the entry of $L$ in row $i$ column $j$. 
Set
\dsp{L_{i,j}=\sum_{k} L_{i,j,k}\xi^k}.
By equation \eqref{negpart},   $(L)_-$ has order at most $0$.   (This 
is somewhat
surprising.)  Thus $L_{1,2}$ has order at most $c_1-c_2$, so 
$L_{1,2,c_1-c_2} $ is a section of $\O$.    Then
\begin{equation}\label{flow5}\partial_m(L_h)=[M,L_h]\ 
,\end{equation}
where $M$ is the following modification of $(L)_+$:
\begin{equation}M=(L)_{+}-\begin{pmatrix}0&0\\0&L_{1,2,c_1-c_2}\end{pmatrix}\ .\end{equation} 
\end{Ex}

\begin{Ex}     Assume now that ${\vec Q}=(Q_1,Q_2,Q_3)$,  with $Q_1\in S$, 
and $\psi$ partially normalized
by
\dsp{\tilde\psi|_{\vec Q}=\begin{pmatrix}1&a&b\\0&c&d\\0&0&e\end{pmatrix}}.   Then
\dsp{R=\begin{pmatrix}0&\frac{a'}c&-\frac{a'd}{ec}+\frac{b'}{e}\\
0&\frac{c'}c&-\frac{c'd}{ec}+\frac{d'}{e}\\
0&0&\frac{e'}e\end{pmatrix}}.

If $a$ is invertible,  take $c=a$,   yielding
\dsp{\begin{pmatrix}0&\frac{a'}a&-\frac{a'd}{ea}+\frac{b'}{e}\\
0&\frac{a'}a&-\frac{a'd}{ea}+\frac{d'}{e}\\
0&0&\frac{e'}e\end{pmatrix}}.

Then if $d$ is invertible,   take $e=d$. This gives equation \eqref{flow5}, with
\begin{equation}
M=(L)_+-\begin{pmatrix}0&0&0\\0&L_{1,2,c_1-c_2}&0\\0&0&L_{2,3,c_2-c_3}+
L_{1,2,c_1-c_2}\end{pmatrix}\end{equation}

On the other hand,  if $d-b$ is invertible,   we may also take 
$e=d-b$. This gives 
\begin{equation}
M=(L)_+-\begin{pmatrix}0&0&0\\0&L_{1,2,c_1-c_2}&0\\0&0&
L_{2,3,c_2-c_3}-L_{1,3,c_1-c_3}\end{pmatrix}\end{equation}

\end{Ex}

The general procedure starts with a good projective basis $\psi$,  
defined over some open subset $\Cal U\subset Y$, and with $S$ 
consisting of a single point.
The procedure recursively normalizes the rows, possibly shrinking $\Cal U$, and 
introducing more elements into
$S$ only as a last resort.
\bigskip

\begin{Proc}[ {\bf Normalize($\psi$)}]\label{normalize}
\

\medskip

{\bf Initialization:}   Set $S=\{Q_1\}$, so that the line bundle  
$\Cal L_1$ is the trivial.   Then   normalize the first row of  $\tilde\psi|_{\vec Q}$ to
be 
$(1,*,...,*)$.
\medskip

{\bf Loop:}   For $\ell\ge 1$, assume the first
$\ell$ elements of $\psi$ have been normalized.   Let
$\psi^{(\ell)}$ denote the upper left-hand $\ell\times\ell$ 
block.
Set
  $R^{(\ell)}={(\widetilde{\psi^{(\ell)}}|_{\vec Q})}'{(\widetilde{\psi^{(\ell)}}|_{\vec Q})}^{-1}$, 
where $(\cdot)'$ denotes $\nabla_{\xi}(\cdot)$.  Set
\begin{equation}
\psi^{(\ell+1)}|_{\vec Q}=\begin{pmatrix}\psi^{(\ell)}|_{\vec Q}&v\\0&x\end{pmatrix}\ ,
\end{equation}
where $x$ is a section of $\Cal L_{\ell+1}$,   yet to be determined.
Then
\begin{equation}
R^{(\ell+1)}=\begin{pmatrix}R^{(\ell)}&-R^{(\ell)}\frac 
vx+\frac{v'}x\\0&\frac{x'}x\end{pmatrix}\ .
\end{equation}
Let $j$ be an integer such that \dsp{v=\begin{pmatrix}v_1\\0\\ 
\vdots\\ 0\end{pmatrix}},   where $v_1$ is a vector of length $j$.

If possible,  choose  a  vector of constants $c=(c_1,...,c_j)$ such 
that   there exists
$k\in \Gamma(\Cal U,\O)$ such that
$
cR^{(j)}=kc
$.
Let
$n$ denote $c$ padded at the end with $l-j$ zeroes.     Then 
$nR^{(\ell)}=(kc,*)=kn+m$,  where $m\cdot v=0$.

  If $n\cdot v$ is not the $0$ section of  $\Cal L_{\ell+1} $, shrink 
$\Cal U$ so that  $n\cdot v$ is invertible.  Declare $x=n\cdot 
v=c\cdot v_1$.   Note that
$x'=n\cdot v'$. Then

\begin{align}
n\cdot(-R^{(\ell)}\frac vx+\frac{v'}x)
&=-k+\frac{x'}x
\end{align}

Thus,  for some $i\le j$,
\begin{equation}
{R^{(\ell+1)}}_{\ell+1,\ell+1}=\sum_{m=1}^j c_i 
{R^{(\ell+1)}}_{m,\ell+1} +{R^{(\ell+1)}}_{i,i}
\end{equation}

If no suitable vector $c$ can be found,  adjoin $Q_{\ell+1}$ to $S$. 
This makes  $\Cal L_{\ell+1}$ the trivial bundle.    Then declare 
$x=1$.
{\bf end}

\end{Proc}
It is clear that  procedure \ref{normalize} solves problem \ref{nproblem3}.

\subsubsection{Formal setting}   Procedure \ref{normalize} may be applied in
the formal setting in 
which expressions that are not formally zero are presumed to be 
nonzero sections.   

\begin{Prop}\label{formalsetting}
In the formal setting it is
possible to   follow  procedure \ref{normalize} in such a way that at the end, the 
number of points
in
$S$ is equal to the multiplicity of $c_1$.   \end{Prop}
I leave the proof of proposition \ref{formalsetting} to the reader. One can get 
the
gist of it  by considering the  case that each $c_i$
occurs with mulitiplicity one.    Starting with $\psi$ such that
\begin{equation}
\tilde\psi|_{\vec Q}=\begin{pmatrix}1&a_{1,2}&&...&\\0&a_{2,2}&a_{2,3}&...&\\&...&&&a_{d-1,d}\\&&&&a_{d,d}\end{pmatrix}\ 
,\end{equation}
one possibility is to recursively normalize $\psi$ so that 
$a_{k,k}=a_{k-1,k} $.   Another possibility is
$a_{k,k}=a_{1,k}-\sum_{j=2}^{k-1} a_{j,k} $.   
There are intermediate choices as well,  so the possiblilities grow 
exponentially with $d$.   For instance,   with $d=4$,   one could 
choose
\begin{align}a_{2,2}&=a_{1,2}\notag\\
a_{3,3}&=a_{2,3}\notag\\
a_{4,4}&=a_{2,4}-a_{3,4}\end{align}
or
\begin{align}a_{2,2}&=a_{1,2}\notag\\
a_{3,3}&=a_{1,3}-a_{2,3}\notag\\
a_{4,4}&=a_{3,4}\ .\end{align}

By virtue of 
proposition \ref{formalsetting},  one may expect to be able to reduce $S$ to a single point
in the event that $c_1$ occurs with multiplicity
one.

\section{Examples in which one may reduce $S$ to a point}\label{sec:examples}   

\begin{Ex}[{\bf Fano Surface}]  Let $F$ be the family of lines on
a  smooth cubic threefold.  (See \cite{B,CG} for background.) 
Then $F$  is a smooth surface.    For
a general line
$\ell$ on the threefold,   let
$D$ be the family of lines on the threefold meeting   $\ell$ 
transversely.   Then $D$ is a smooth, ample curve, of genus 11.   Let 
$W_0=H^0(D,\Cal N)$.
Then $W_0$ is canonically isomorphic to the tangent space 
$T_{\ell}F$.   In particular,   $W_0$ is two-dimensional and 
basepoint free.   The degree of $\Cal N$
is 5.
Given $\ell'\in D$,   the plane spanned by $\ell$ and $\ell'$ meets 
the threefold in $\ell$,  $\ell'$ and a third
line, which we may denote by $\sigma(\ell')$.
Then $\sigma$ is an involution with no fixed points,   and
the natural map $\Pic^0(F)\to\Pic^0(D)
$ imbeds
$\Pic^0(F)$ as the Prym variety of $(D,\sigma)$.  For a general 
degree-zero line bundle $\Cal
L\in\Pic^0(F)$ one has the following invariants:
\begin{equation}\begin{matrix}h^0(\Cal L(kD))=h^1(\Cal L(kD))=0\ ,&k<2,\\\\
h^1(\Cal L(kD))=h^2(\Cal L(kD))=0\ ,&k\ge 2.
\end{matrix}\end{equation}
By theorem \ref{Locally Free at the generic point},  for any $S$ and 
any space $W\subset H^0(F,\O(D)/\Cal I_S)$ mapping isomorphically 
onto $W_0$,
$\ft(\O(*D))$ is a locally free $\B_W$-module in a neighborhood of 
$\Cal L$.     Furthermore,
\begin{equation}h^0(\Cal L(2D))=1\ , h^0(\Cal L(3D))=6\ , h^0(\Cal L(4D))=16\ ,
\end{equation}
from which it follows that $c_1=2$,   $c_2=c_3=c_4=3$ and $c_5=4$. 
In particular,   $c_1$ occurs with multiplicity one.    Thus, if 
$\psi$ is  any    good
basis,   $\tilde\psi|_{\vec Q}$ has the following form:
\begin{equation}
\tilde\psi|_{\vec Q}=\begin{pmatrix}\alpha_{1,1}&\alpha_{1,2}&\alpha_{1,3}&\alpha_{1,4}&\alpha_{1,5}\\0&\alpha_{2,2}&0&0&\alpha_{2,5}\\
0&0&\alpha_{3,3}&0&\alpha_{3,5}
\\
0&0&0&\alpha_{4,4}&\alpha_{4,5}\\
0&0&0&0&\alpha_{5,5}\end{pmatrix}\end{equation}
Let $L=(L_g)^{-1}L_f $ and define $L_{i,j,k}$ as in example \ref{flowexample}.
We may set $S=\{Q_1\}$ and set $\alpha_{1,1}=1$.
By choosing $\alpha_{2,2}=\alpha_{1,2}$,
$\alpha_{3,3}=\alpha_{1,3}$,
$\alpha_{4,4}=\alpha_{1,4}$ and
$\alpha_{5,5}=\alpha_{4,5}$, we find that the system
$$\partial_m(L_h)=[M,L_h]$$ 
 may be 
solved by rational functions on $\Pic^0(F)$,
where 
\begin{equation}\kern-20ptM=(L)_{+}-\begin{pmatrix}0&0&0&0&0\\0&L_{1,2,-1}&0&0&0\\
0&0&L_{1,3,-1}&0&0\\
0&0&0&L_{1,4,-1}&0\\
0&0&0&0&L_{4,5,-1}+L_{1,4,-1}\end{pmatrix}\end{equation}
Another possibility is $\alpha_{2,2}=\alpha_{1,2}$,
$\alpha_{3,3}=\alpha_{1,3}$,
$\alpha_{4,4}=\alpha_{1,4}$ and
$\alpha_{5,5}=\alpha_{1,5}-\alpha_{2,5}-\alpha_{3,5}-\alpha_{4,5}$. 
We then get 
\begin{align}\label{flow7}&M=\\&(L)_{+}-\begin{pmatrix}0&0&0&0&0\\0&L_{1,2,-1}&0&0&0\\
0&0&L_{1,3,-1}&0&0\\
0&0&0&L_{1,4,-1}&0\\
0&0&0&0&L_{1,5,-2}-L_{2,5,-1}-L_{3,5,-1}-L_{4,5,-1}\end{pmatrix}\ 
.\notag\end{align}
Both possibilities are allowed for generic threefold and generic line
on it.\end{Ex}

\begin{Ex}[{\bf Ruled Surface}] 
Let  $C$ be a smooth hyperelliptic curve, \dsp{C\morph{\sigma}\BBB P^1}, with
$g=g(C)\ge 3$.  
Consider an extension
\begin{equation}
  \label{eq:extension}
  0\to\O\morph t\Cal E\to\Cal L\to 0\ .
\end{equation}
where $\Cal L=\sigma^*(\O_{\BBB P^1}(1))$.   
Set \dsp{X=\BBB P(\Cal E)\morph{\pi}C}.

If the  extension class of $\Cal E$ is chosen generally,   then 
the two-dimensional vector space
$H^0(C,\Cal L)$ will inject into
$H^1(C,\O)$.   Furthermore,  the relative $\O(1)$ on $X$ will be ample. 
Then,  if we let $Z\subset X$ be the curve defined by
$\{t=0\}$,  $Z$ is ample.   Denoting by $\Cal N$
the normal bundle of $Z$,   one recovers the sequence \eqref{eq:extension} by applying
$\pi_*$ to the sequence
\begin{equation}
    \label{eq:normalsequence}
 0\to\O\to\Cal \O(Z)\to\Cal N\to \O\  .  
  \end{equation}
In particular,  $\dim H^0(Z,\Cal N)=2$, and one has an injection
\begin{equation}
0\to H^0(Z,\Cal N)\to H^1(X,\O)=H^1(C,\O)\ .\end{equation} 

\begin{Prop}   Let $\Cal M$ be a general line bundle on $C$, of degree
$g-1\le d\le g+1$.  Let $b$ be an integer.  Then
\begin{equation}
    \label{eq:vanishing}
    H^i(X,\pi^*(\Cal M)(bZ))=0
  \end{equation}
in any of the following circumstances:

  1.   $i=1$,

2. $b\ge -1$ and $i=2$,

3. $b\le -1$ and $i=0$.

\end{Prop}
\begin{pf}

Case $i=1$:

When $b=0$,   we are asserting that $H^1(C,\Cal M)=0$,
which holds generically for the stated degrees.
For $b>0$,  we have $H^1(X,\Cal N^b\otimes\pi^*(\Cal M))=H^1(C,\Cal L^b\otimes\Cal M)=0$.
Then statement 1 follows for all positive $b$ by induction,  using the
exact sequence
\begin{equation}
    \label{eq:vanishbyinduction}
 0\to\O(b-1)\otimes\pi^*\Cal M\to\Cal \O(b)\otimes\pi^*\Cal M
\to\Cal N^b\otimes\pi^*\Cal M\to \O\ .   
  \end{equation}
If $b\le -1$, then
$deg(\Cal L^b\otimes\Cal M)\le g-1$. Then 
use \eqref{eq:vanishbyinduction} again, invoking this time the fact
that $H^0(X,\Cal N^b\otimes\pi^*\Cal M)=0$.   This finishes the case $i=1$.

The other two cases follow from the fact that for $b\ge -1$, $R^1\pi_*(\O(b))=0$, while for
$b\le -1$, $R^0\pi_*(\O(b))=0$.

\end{pf}

Let  $g-1\le d\le g+1$ and let $\kappa=d-(g-1)$.
Let $\chi(b)$ denote the Euler characteristic of  $\pi^*\Cal M(b)$,  
for $deg(\Cal M)=d$.   Then one easily computes that
\begin{equation}
  \label{eq:eulercharacteristic}
  \chi(b)=\kappa+(\kappa+1)b+b^2\ .
\end{equation}
We therefore get the following theorem.
\begin{Thm}
  \label{punchline}
The hypotheses of theorem 3.2 hold with $X$ and $Z$ as above,  
 taking $\Cal F$ to be the pullback from $C$ of any line 
bundle $\Cal M$ with $\kappa=0,1,2$, where $\kappa=deg(\Cal M)-(g-1)$.   The filtered $\B_W$-module
obtained in this way has two generators of the following orders:

0.   both of order $1$ if $\kappa=0$,

1.   one of order $0$ and one of order $1$ if $\kappa=1$,

2.  both of order $0$ if $\kappa=2$.
\end{Thm}

With $\kappa=1$, i.e.,  $deg(\Cal M)=g$,  one may reduce $S$ to a single point.   Then one obtains
the flows 
$$\partial_m(L_h)=[M,L_h]\ ,$$ 
where,  $M$ is as  in example
\ref
{flowexample},  with $c_1-c_2=-1$:
\begin{equation}M=(L)_{+}-\begin{pmatrix}0&0\\0&L_{1,2,-1}\end{pmatrix}\ .\end{equation}

\end{Ex}

\bigskip
{\sc
Department of Mathematics, University of Georgia, Athens, GA
30602}

{\it E-mail address:}
rothstei@@math.uga.edu

\end{document}